\documentclass[fleqn,twoside,twocolumn,nofootinbib,showkeys]{revtex4} 
\usepackage[nocpr]{ujp} 

\newcommand{\es}{\end{multicols}}
\newcommand{\ed}{\end{document}}
\newcommand{\be}{\begin{equation}}
\newcommand{\ee}{\end{equation}}
\newcommand{\bc}{\begin{center}}
\newcommand{\ec}{\end{center}}
\newcommand{\ba}{\begin{array}}
\newcommand{\ea}{\end{array}}

\begin{document}
\title[HOMFLY Polynomial Invariants of Torus Knots]
{HOMFLY POLYNOMIAL INVARIANTS\\ OF TORUS KNOTS AND BOSONIC \boldmath$(q,p)$-CALCULUS$^1$}%
\author{A.M.~Pavlyuk}
\affiliation{Bogolyubov Institute for Theoretical Physics, Nat. Acad. of Sci. of Ukraine}
\address{14b, Metrolohichna Str., Kyiv 03143, Ukraine}
\email{pavlyuk@bitp.kiev.ua}

\udk{539.12; 517.984}  \razd{\secix}

\autorcol{A.M.\hspace*{0.7mm}Pavlyuk}%

\setcounter{page}{1178}%

\begin{abstract}
For the one-parameter Alexander (Jones) skein relation we introduce
the Alexander (Jones) ``bosonic'' $q$-numbers, and for the
two-parameter HOMFLY skein relation we propose the \mbox{HOMFLY}
 ``bosonic'' $(q,p)$-numbers (``bosonic'' numbers connected with
deformed bosonic oscillators). With the help of these  deformed
 ``bosonic'' numbers, the corresponding skein relations can be
reproduced. Analyzing the introduced ``bosonic'' numbers, we point
out two ways of obtaining the two-parameter HOMFLY skein relation
(``bosonic'' $(q,p)$-numbers) from the one-parameter Alexander and
Jones skein relations (from the corresponding ``bosonic''
$q$-numbers). These two ways of obtaining the HOMFLY skein relation
are equivalent.

\end{abstract}


\maketitle

\section{Introduction}

The knot theory is substantially based on the axioms of skein
relation and normalization~\cite{Co} allowing one to describe every
knot and link by a definite polynomial. These polynomials form the
set of polynomial invariants. The goal of this paper is to show that
every of the three polynomial invariants (Alexander, Jones, HOMFLY)
can be put into correspondence to the definite ``bosonic''
$q$-numbers/$(q,p)$-numbers $([n]^{A},\ [n]^{V},\ [n]^{H})$, which
allow one to reproduce the corresponding skein relation. Comparing
these deformed numbers, we find the rule of obtaining the
two-parameter HOMFLY polynomial invariants of knots and links from
one-parameter (Alexander, Jones) polynomial invariants.

The so-called $(q,p)$-numbers, which are an important ingredient
of $(q,p)$-calculus, generalizing the well-known
$q$-calculus~\cite{Ka}, appear in the connection with
$(q,p)$-deformed bosonic oscillators~\cite{CJ}.
 Recently, with the help of these $(q,p)$-numbers, we have shown for
 the important case of torus knots that the generalized (two-variable)
 Alexander polynomials can be obtained~\cite{GP1, GP2}, as well as the generalized
 (two-variable) Jones polynomials~\cite{Pa-agg, Pa}.

\section{Skein Relations}

The Alexander polynomials $\Delta(t)$ for knots and links are
defined by the Alexander skein relation~\cite{Al} together with the
normalization condition for \mbox{the unknot:}\footnotetext[1]{This
work is the contribution to Proceedings of the International
Conference ``Quantum Groups and Quantum Integrable Systems''.}
  \be\label{alex-skein}
\Delta_{+}(t)-\Delta_{-}(t)=(t^{1\over2}-t^{-{1\over2}})\Delta_{O}(t),
\Delta_{\rm unknot}{=}1.
 \ee
The Jones polynomials $V(t)$ are described by the following skein
relation and normalization condition~\cite{Jo}:
\be\label{jones-skein}
t^{-1}V_{+}(t)-tV_{-}(t)=(t^{1\over2}-t^{-{1\over2}})V_{O}(t),
V_{\rm unknot}{=}1.
 \ee
The HOMFLY polynomials $H(a, z)$ are introduced in the following way~\cite{FY}:
\[
 a^{-1}H_{+}(a, z)-aH_{-}(a, z)=zH_{O}(a, z),\ H_{\rm unknot}{=}1.
 \]
For our goal, it is necessary to make a change of the variable
$z=t^{1\over2}-t^{-{1\over2}}.$ Thus, the HOMFLY skein relation
can be rewritten in the form
\be\label{homfly-skein}
\ a^{-1}H_{+}(a, t)-aH_{-}(a, t)=(t^{1\over2}-t^{-{1\over2}})H_{O}(a, t).
 \ee

Let us write the skein relation in the general form
 \be\label{skein}
P_{L_{+}}(t)=l_{1}P_{L_{O}}(t)+l_{2}P_{L_{-}}(t),\ee where $l_{1}$
and $l_{2}$ are coefficients. The capital letter ``L'' stands for
``Link'' and denotes one of the two: knot or link (unknot belongs
to knots). Here, three polynomials $P_{L_{+}}(t)$, $P_{L_{O}}(t),$
$P_{L_{-}}(t)$ correspond to the overcrossing Link $L_{+}$
(``overcrossing'' refers to a chosen crossing of the Link), zero
crossing Link $L_{O},$ and undercrossing Link $L_{-}.$ Thus,
applying the surgery operation of elimination of a crossing to an
initial Link $L_{+},$ one obtains a simpler Link $L_{O}$. The Link
$L_{-}$ is obtained from the same initial Link  $L_{+}$ by the
another surgery operation of switching of the crossing.

Consider the simplest torus knots $T(2m+1, 2)$ and torus links
$L(2m, 2),$ where $m=0,1,2,3,\mbox{...}\,$. The common notation
for these torus knots and links  $L_{n, 2}$ corresponds to torus
knots, if $n$ is odd, and to torus links for even $n.$ The surgery
operation of elimination turns  $L_{n, 2}$  into $L_{n-1, 2}$, and
the switching operation turns it into $L_{n-2, 2}$. Because of it,
the very important property follows from~(\ref{skein}), namely,
the series of polynomials $P_{L_{n, 2}}(t)$ is characterized by
the recurrence relation
\be\label{skein1} P_{L_{n+1, 2}}(t)=l_{1}P_{L_{n,
2}}(t)+l_{2}P_{L_{n-1, 2}}(t),\ee which repeats itself in the
skein relation~(\ref{skein}). We now rewrite
formula~(\ref{skein1}) in a simpler notation
 \be\label{skein2}
P_{n+1, 2}(t)=l_{1}P_{n, 2}(t)+l_{2}P_{n-1, 2}(t).
 \ee
For the polynomials $P_{L_{n, 2}}(t)$ with odd $n$, 
relation~(\ref{skein2}) yields the recurrence relation referred
only to the torus knots $T(2m+1, 2)$
 \be\label{skein3}
P_{n+2, 2}(t)=k_{1}P_{n, 2}(t)+k_{2}P_{n-2, 2}(t),
 \ee
where the coefficients $k_{1}$ and $k_{2}$ are expressed through
$l_{1}$ and $l_{2}$~\cite{GP1} as
\be\label{k}
k_{1}=l_{1}^{2}+2l_{2},\quad k_{2}=-l_{2}^{2}.\ee
From~(\ref{skein3}) and the normalization condition
\be\label{norm} P_{1, 2}=1 \ee we obtain for the trefoil
\be\label{32} P_{3,2}=k_{1}+k_{2}. \ee

\section{``Bosonic'' \boldmath$(q,p)$-numbers}

The one-parameter ``bosonic'' $q$-number (structural function)
characteristic of a Biedenharn--Macfarlane deformed bosonic
oscillator corresponding to an integer $n$
 is defined as~\cite{Bi,Ma}
 \be \label{q-def} [n]_{q}={{\textstyle
{q^{n}-q^{-n}}}\over{\textstyle {q-q^{-1}}}},
 \ee
  with $q$ to be  a parameter.
Some of the $q$-numbers are
\[
[1]_{q}=1,\quad  [2]_{q}=q+q^{-1},\]
 \[ [3]_{q}=q^{2}+1+q^{-2}, \quad  [4]_{q}=q^{3}+q+q^{-1}+q^{-3}, \mbox{...}\,. \]
  The recurrence relation for~(\ref{q-def})
looks as
\be\label{qp-rec} [n+1]_{q}=(q+q^{-1})[n]_{q}-[n-1]_{q}.\ee

The two-parameter ``bosonic'' $(q,p)$-number corresponding to the
integer number $n$ is defined as~\cite{CJ}
 \be \label{qp-def}
[n]_{q,p}={{\textstyle {q^{n}-p^{n}}}\over{\textstyle {q-p}}}, \ee
where $q,p$ are  parameters. If $p=q^{-1},$ then $[n]_{q,p}=$
$=[n]_{q}.$ Some of the $q,p$-numbers are given below:
 \[ [1]_{q,p}=1,\quad [2]_{q,p}=q+p,
  \]\vspace*{-7mm}
 \[ [3]_{q,p}{=}q^{2}+qp+p^{2},\quad   [4]_{q,p}{=}q^{3}+q^{2}p+qp^{2}+p^{3}, \mbox{...}\,.\]
The recurrence relation for $q,p$-numbers is
\be\label{qp-rec} [n+1]_{q,p}=(q+p)[n]_{q,p}-qp[n-1]_{q,p}.\ee

\section{Alexander ``Bosonic'' \boldmath$q$-numbers: \boldmath$ [n]^{A}$ }

From the Alexander skein relation~(\ref{alex-skein}) in the form~(\ref{skein})
 \be\label{alex-skein2}
\Delta_{+}(t)=(t^{1\over2}-t^{-{1\over2}})\Delta_{O}(t)+\Delta_{-}(t),
 \ee
  one has the ``Link coefficients''
 \be\label{alex-l}
l_{1}^{A}=t^{1\over2}-t^{-{1\over2}},\quad  l_{2}^{A}=1.
 \ee
From~(\ref{alex-skein2}), we have the recurrence relation for
Alexander polynomials of torus knots and links $L_{n, 2}$ (by
analogy to~(\ref{skein2}))
\be\label{alex-skein3}
\Delta_{n+1,2}(t)=(t^{1\over2}-t^{-{1\over2}})\Delta_{n,2}(t)+\Delta_{n-1,2}(t).
 \ee
Using~(\ref{k}) and~(\ref{alex-l}), we obtain the ``knot
coefficients''
\be\label{alex-k} k_{1}^{A}=t+t^{-{1}},\quad k_{2}^{A}=-1.
 \ee
Therefore, the recurrence relation for Alexander polynomials of
torus knots $T(2m+1, 2)$ looks as
\be\label{alex-skein4}
\Delta_{n+2,2}(t)=(t+t^{-{1}})\Delta_{n,2}(t)-\Delta_{n-2,2}(t).
 \ee
Comparing~(\ref{alex-skein4}) and~(\ref{qp-rec}) allows us to put
(what we call) the Alexander ``bosonic'' $q$-numbers $ [n]^{A}$ into
correspondence to~(\ref{alex-skein4}). Indeed, from $
q+p=t+t^{-{1}}$, $qp=1,$ we have $ q=t$, $p=t^{-{1}}.$ Therefore,
relation~(\ref{qp-def}) yields
\be \label{alex-b} [n]^{A}={{\textstyle
{t^{n}-t^{-{n}}}}\over{\textstyle {t-t^{-{1}}}}},\quad t\equiv q, \ee
which
coincides with $q$-numbers of Biedenharn and Macfarlane.

\section{Jones ``Bosonic'' \boldmath$q$-numbers: \boldmath$ [n]^{V}$ }

From the Jones skein relation~(\ref{jones-skein}) in the form~(\ref{skein})
 \be\label{jones-skein2}
V_{+}(t)=t(t^{1\over2}-t^{-{1\over2}})V_{O}(t)+t^{2}V_{-}(t),
 \ee
we have the ``Link coefficients''
  \be\label{jones-l}
l_{1}^{V}=t^{3\over2}-t^{{1\over2}},\quad  l_{2}^{V}=t^{2},
 \ee
and, correspondingly, find the ``knot coefficients''
\be\label{jones-k} k_{1}^{V}=t^{3}+t,\quad  k_{2}^{V}=-t^{4}.
 \ee
The recurrence relation for the Jones polynomials of torus knots
$T(2m+1, 2)$ has the form
\be\label{jones-skein4}
V_{n+2,2}(t)=(t^{3}+t)V_{n,2}(t)-t^{4}V_{n-2,2}(t).
 \ee
Comparing~(\ref{jones-skein4}) and~(\ref{qp-rec}), we obtain what we
call the Jones ``bosonic'' $q$-numbers
\be \label{jones-b} [n]^{V}={{\textstyle
{t^{3n}-t^{{n}}}}\over{\textstyle {t^{3}-t}}}, \quad t\equiv q, \ee

\section{HOMFLY ``Bosonic'' \boldmath$(q,p)$-numbers: \boldmath$ [n]^{H}$ }

The HOMFLY skein relation~(\ref{homfly-skein}) in the form~(\ref{skein})
 \be\label{homfly-skein2}
\ H_{+}(a, t)=a(t^{1\over2}-t^{-{1\over2}})H_{O}(a, t)+a^{2}H_{-}(a, t)
 \ee
 gives the ``Link coefficients''
 \be\label{homfly-l}
l_{1}^{H}=a(t^{1\over2}-t^{-{1\over2}}),\quad  l_{2}^{H}=a^{2}.
 \ee
From whence, we find the ``knot coefficients''
\be\label{homfly-k} k_{1}^{H}=a^{2}(t+t^{-{1}}),\quad  k_{2}^{H
}=-a^{4},\ee which are used to introduce the HOMFLY ``bosonic''
$(q,p)$-numbers  according to the relation
\be\label{n-bv}
[n+1]^{H}=k_{1}^{H}[n]^{H}+k_{2}^{H}[n-1]^{H}. \ee
 Comparing~(\ref{n-bv}) and~(\ref{qp-rec}),
from $ q+p=k_{1}^{H}$, $qp=$ $=-k_{2}^{H},$ we have $q=at$ and
$p=at^{-{1}}.$  It follows from~(\ref{qp-def}) what we call the
HOMFLY ``bosonic'' $(q,p)$-numbers
\be \label{homfly-b} [n]^{H}=a^{2(n-1)}{{\textstyle
{t^{n}-t^{-{n}}}}\over{\textstyle {t-t^{-{1}}}}},\quad t\equiv q,\ a\equiv p. \ee

\section{Alexander Skein Relation\\ from Alexander ``Bosonic'' \boldmath$q$-numbers}

In Section 4, we obtained the Alexander
 ``bosonic'' $q$-numbers from the Alexander skein relation.
In this section, moving in the opposite direction, we obtain the
Alexander skein relation~(\ref{alex-skein2}) from  the Alexander
 ``bosonic'' $q$-numbers $[n]^{A}$~(\ref{alex-b}). First,
comparing~(\ref{qp-def}) and~(\ref{alex-b}), we find $q=t$,
$p=t^{-1}$. Putting it into~(\ref{qp-rec}), one has the recurrence
relation for the Alexander ``bosonic'' $q$-numbers:
\be\label{n-alex}
[n+1]^{A}=(t+t^{-1})[n]^{A}-[n-1]^{A}. \ee%
From whence, we have ``knot coefficients''~(\ref{alex-k}):
\[
k_{1}^{A}=t+t^{-{1}},\quad k_{2}^{A}=-1.\] According to~(\ref{k}),
the ``Link coefficients'' are
\be\label{l-k} l_{2}=+(-k_{2})^{1\over2},\quad   l_{1}=
+(k_{1}-2l_{2})^{1\over2}.\ee Thus, we obtain
$l_{2}^{A}$ and $l_{1}^{A},$ which coincide with~(\ref{alex-l}).
Putting them into~(\ref{skein}) leads to the Alexander skein
relation~(\ref{alex-skein2}).

In a similar way,  the Jones ``bosonic'' $q$-numbers $
[n]^{V}$~(\ref{jones-b}) yield the Jones skein
relation~(\ref{jones-skein2}), and the HOMFLY skein
relation~(\ref{homfly-skein2}) follows from the HOMFLY ``bosonic''
$(q,p)$-numbers $ [n]^{H}$~(\ref{homfly-b}).

\section{HOMFLY invariants\\ from Alexander and Jones invariants}

In this section, we consider how to build two-parameter HOMFLY
polynomial invariants on the basis of one-parameter Alexander or
Jones ones. To formulate proper rule, we compare the HOMFLY
 ``bosonic'' $(q,p)$-numbers $ [n]^{H}$~(\ref{homfly-b}) and the
Alexander ``bosonic'' $q$-numbers $ [n]^{A}$~(\ref{alex-b}):
\be\label{homfly-alex} [n]^{H}=a^{2(n-1)}[n]^{A}. \ee
Then, the first way of obtaining the HOMFLY $(q,p)$-numbers reduces
 to introducing the second variable in the form of a
multiplier $a^{2(n-1)}$ before $[n]^{A}$. In  the case of the Jones
 ``bosonic'' $q$-numbers $[n]^{V},$ the multiplier looks as
$(aq)^{2(n-1)}$:
\be\label{homfly-jones}
[n]^{H}=(aq)^{2(n-1)}[n]^{V}.
\ee

We suggest another way of obtaining the HOMFLY skein relation with
the help of the ``bosonic'' $(q,p)$-numbers $ [n]^{H_{1}}$. To this
end, we make the substitution in~(\ref{alex-b}):\vspace*{-1mm}
\[t\rightarrow q,\quad t^{-1}\rightarrow p^{-1}\]
 and, thus,
\be \label{homfly-b1} [n]^{H_{1}}={{\textstyle
{q^{n}-p^{-{n}}}}\over{\textstyle {q-p^{-{1}}}}}.
 \ee
One more substitution
\[q^{1\over4}p^{-{1\over4}}\rightarrow a,\quad q^{1\over2}p^{{1\over2}}\rightarrow t,\] in~(\ref{homfly-b1})
turns $[n]^{H_{1}}$ into $[n]^{H},$ which proves their
equi\-valence.

In  the case of the Jones ``bosonic'' $q$-numbers $ [n]^{V},$ the
substitution
\[t^{3}\rightarrow q^{3},\quad t\rightarrow p\]
in~(\ref{alex-b}) turns it into
\be \label{homfly-b2} [n]^{H_{2}}={{\textstyle
{q^{3n}-p^{{n}}}}\over{\textstyle {q^{3}-p}}}. \ee By substituting
\[q^{3}p\rightarrow a^{4},\quad q^{3\over2}p^{-{1\over2}}\rightarrow
t,\] in~(\ref{homfly-b2}), we turn $[n]^{H_{2}}$ into
$[n]^{H}.$\vspace*{2mm}

\section{Concluding Remarks}

The introduced Alexander and Jones ``bosonic'' $q$-numbers and the
HOMFLY ``bosonic'' $(q,p)$-numbers give possibility to operate with
these deformed numbers instead of operating with the corresponding
skein relations, which is believed to be much easier. We also hope
that the dealing with the deformed numbers instead of the skein
relations will promote the finding of new polynomial invariants of
knots and links. It should be mentioned that the problem of
searching for the Reidemeister moves in terms of $(q,p)$-calculus
arises.

\vskip3mm {\it This work was partially supported by the Special
Programme of the Division of Physics and Astronomy of the NAS of
Ukraine.}



\vspace*{-5mm}
\rezume{%
A.M.\,Павлюк} {ПОЛІНОМІАЛЬНІ ІНВАРІАНТИ\\ ХОМФЛІ ДЛЯ ТОРИЧНИХ
ВУЗЛІВ\\
І БОЗОННЕ $(q,p)$-ЧИСЛЕННЯ} {Для однопараметричного
скейн-співвідношення Алек\-сан\-де\-ра (Джонса) введено ``бозонні''
$q$-числа Алек\-сандера (Джон\-са),  а для двопараметричного
скейн-співвідношення Хомфлі~-- ``бозонні'' $(q,p)$-числа Хомфлі
(``бозонні'' числа пов'язані з деформованими бозонними
осциляторами). За допомогою цих  деформованих ``бозонних'' чисел
можна відновити відповідні скейн-співвідношення. Аналізуючи введені
``бозонні'' числа, ми вказуємо на два способи отримання
двопараметричного скейн-спів\-від\-но\-шен\-ня Хомф\-лі
(``бозонних'' $(q,p)$-чисел) із однопараметричних
скейн-спів\-від\-но\-шень Александера і Джонса (із відповідних
``бозонних'' $q$-чисел). Ці два способи отримання
скейн-співвідношення Хомфлі \mbox{еквівалентні.}}

\end{document}